\documentclass[11pt]{article}

\usepackage{amsmath}
\usepackage{amssymb}
\usepackage{indentfirst}
\usepackage{graphics} 
\usepackage{color}

\makeatletter

\@addtoreset{equation}{section}
\makeatother
\def\<{\langle}
\def\>{\rangle}

\newtheorem{lem}{Lemma}[section]
\newtheorem{theo}{Theorem}[section]
\newtheorem{rem}{Remark}[section]
\newtheorem{pro}{Proposition}[section]

\makeatletter
   
   \@addtoreset{equation}{section}
\makeatother

\begin{document}
\title{\bf A Note on Decay Rates of the Local Energy \\ for Wave Equations with Lipschitz Wavespeeds}
\author{Ruy Coimbra Char\~ao \\{\small Department of Mathematics} \\{\small Federal University of Santa Catarina} \\ {\small 88040-270, Florianopolis, Brazil} \\and\\Ryo Ikehata\thanks{Corresponding author: ikehatar@hiroshima-u.ac.jp} \\ {\small Department of Mathematics}\\ {\small Graduate School of Education} \\ {\small Hiroshima University} \\ {\small Higashi-Hiroshima 739-8524, Japan}}
\maketitle
\begin{abstract}
We consider the Cauchy problem for wave equations with variable coefficients in the whole space ${\bf R}^{n}$. We improve the rate of decay of the local energy, which has been recently studied by J. Shapiro \cite{S}, where he derives the log-order decay rates of the local energy under stronger assumptions on the regularity of the initial data.
\end{abstract}
\section{Introduction}
\footnote[0]{Keywords and Phrases: Wave equation; Lipschitz wavespeed; local energy; algebraic decay, non-compactly supported initial data, weighted energy estimate.}
\footnote[0]{2010 Mathematics Subject Classification. Primary 35L05; Secondary 35B40, 35L30.}
We consider in this work  the Cauchy problem associated to the wave equation with variable coefficient  in ${\bf R}^{n}$ ($n \geq 1$) as follow
\begin{equation}
u_{tt}(t,x) -c(x)^{2}\Delta u(t,x) = 0,\ \ \ (t,x)\in (0,\infty)\times {\bf R}^{n},\label{eqn}
\end{equation}
\begin{equation}
u(0,x)= u_{0}(x),\ \ u_{t}(0,x)= u_{1}(x),\ \ \ x\in{\bf R}^{n} ,\label{initial}
\end{equation}
where $(u_{0},u_{1})$ are initial data chosen as
\[
u_{0} \in H^{1}({\bf R}^{n}),\quad u_{1} \in L^{2}({\bf R}^{n}),
\]
and the function $c: {\bf R}^{n} \to {\bf R}$ satisfies the  two assumptions below:\\
{\bf (A-1)}\,$c(x) > 0$ ($x \in {\bf R}^{n}$), $c, c^{-1} \in L^{\infty}({\bf R}^{n})$, $\nabla c \in (L^{\infty}({\bf R}^{n}))^{n}$,\\
{\bf (A-2)}\,there exists a constant $L > 0$ such that $c(x) = 1$ for $\vert x\vert > L$.\\
In particular, the condition (A-1) implies $c \in C^{0,1}({\bf R}^{n})$ (see e.g., \cite[Theorem IX.12]{HB}). Here, we have set
\[u_{t}=\frac{\partial u}{\partial t},\quad  u_{tt}=\frac{\partial^2 u}{\partial t^2}, \quad \Delta = \sum_{j=1}^{n}\frac{\partial^{2}}{\partial x_{j}^{2}}, \quad x = (x_1,\cdots,x_n).\]

Note that solutions and/or functions considered in this paper are all real valued except for some parts concerning the Fourier transform.  

Considering the previous assumptions on the initial data and  $c(x)$  it is known that the problem \eqref{eqn}-\eqref{initial} has a unique weak solution 
\[u \in C([0,\infty);H^{1}({\bf R}^{n}))\cap C^{1}([0,\infty);L^{2}({\bf R}^{n})) =: C_{1}^{n},\]
satisfying the energy conservation property:
\begin{equation}
E_{u}(t) = E_{u}(0),
\end{equation}
where the total energy $E_{u}(t)$ to the equation (1.1) is defined by
\[E_{u}(t) :=\frac{1}{2}\int_{{\bf R}^{n}}\big(\frac{1}{c(x)^{2}}\vert u_{t}(t,x)\vert^{2} + \vert\nabla u(t,x)\vert^{2}\big)dx.\]

Furthermore, the local energy $E_{R}(t)$ on the  zone  $\{\vert x\vert \leq R\}$ ($R > 0$) corresponding to the solution $u(t,x)$ of \eqref{eqn}-\eqref{initial}  is defined by
\[E_{R}(t):=\frac{1}{2}\int_{\vert x\vert \leq R}\big(\frac{1}{c(x)^{2}}\vert u_{t}(t,x)\vert^{2} + \vert\nabla u(t,x)\vert^{2}\big)dx.\]
Also, we set
\[B(x,R) := \{y \in {\bf R}^{n}:\,\vert y-x\vert < R\}.\]
Our main concern of this paper is to obtain a local energy decay estimate with an algebraic decay order. For related important results concerning the local energy decay, one can cite several celebrated papers due to Morawetz \cite{Mora}, Lax-Phillips \cite{LP}, Morawetz-Ralston-Strauss \cite{MRS}, Ralston \cite{R}, Vainberg \cite{V}, and the references therein.

By the way, concerning local energy decay results, quite recently Shapiro \cite{S} announces the following interesting result. It should be mentioned that the decay rate of (1.4) below was first obtained by Burq \cite{B} for smooth perturbations of the Laplacian outside an obstacle. 
\begin{theo}{\rm (Shapiro \cite{S})}\,Let $n \geq 2$, and assume {\rm (A-1)} and {\rm (A-2)}. Suppose that the supports of $u_{0}$ and $u_{1}$ are contained in $B(0,R_{1})$, and $\nabla u_{0} \in (H^{1}({\bf R}^{n}))^{n}$ and $u_{1} \in H^{1}({\bf R}^{n})$. Then for any $R_{2} > 0$, there exists $C > 0$ such that the solution $u$ to {\rm (1.1)}-{\rm (1.2)} satisfies for $t \geq 0$,
\begin{equation}E_{R_{2}}(t) \leq (\frac{C}{\log(2+t)})^{2}\left(\Vert\nabla u_{0}\Vert_{(H^{1}({\bf R}^{n}))^{n}}^{2} + \Vert u_{1}\Vert_{H^{1}({\bf R}^{n})}^{2} \right).\end{equation}
\end{theo} 

Our observation is that Shapiro \cite{S} imposes rather stronger hypothesis on the regularity of the initial data such as\\
{\rm {\bf (I)}}\,the supports of initial data are compact, and as a result $[u_{0},u_{1}] \in H^{2}({\bf R}^{n})\times H^{1}({\bf R}^{n})$.\\
Furthermore, in a sense, \\
{\rm {\bf (II)}}\,the obtained decay order $(\log t)^{-2}$ of the local energy seems to be rather slow.

In this paper, under weaker regularity assumptions on the initial data to modify {\bf (I)}, one obtains faster algebraic decay rate which improves {\bf (II)} in the case when the coefficient $c(x)$ and the parameter $L$ have a special relation. 

Our method is based on the so-called Morawetz identity \cite{Mora}, so we never use the spectral analysis like resolvent estimates. In order to state our results, we introduce the following weighted functional spaces.
\[L^{p,\gamma}({\bf R}^{n}) := \left\{f \in L^{p}({\bf R}^{n}) \; | \; \Vert f\Vert_{p,\gamma} := \left(\int_{{\bf R}^{n}}(1+\vert x\vert^{\gamma})\vert f(x)\vert^{p}dx\right)^{1/p} < +\infty\right\}.\]

Our main results read as follows.
\begin{theo}\,Let $n \geq 3$, and assume {\rm (A-1)} and {\rm (A-2)}. If the initial data $[u_{0},u_{1}] \in H^{1}({\bf R}^{n})\times (L^{2}({\bf R}^{n})\cap L^{1}({\bf R}^{n}))$ further satisfies
\[\int_{{\bf R}^{n}}(1+\vert x\vert)\left(\frac{1}{c(x)^{2}}\vert u_{1}(x)\vert^{2} + \vert\nabla u_{0}(x)\vert^{2}\right)dx < +\infty,\]
then the unique solution $u \in C_{1}^{n}$ to problem {\rm (1.1)}-{\rm (1.2)} satisfies 
\[E_{R}(t) = O(t^{-(1-\eta)})\quad (t \to \infty),\]
for each $R > L$ provided that $\eta := 2L\Vert\frac{1}{c(\cdot)}\Vert_{\infty}\Vert\nabla c\Vert_{\infty} \in [0,1)$.
\end{theo}
\begin{theo}\,Let $n = 2$, and assume {\rm (A-1)} and {\rm (A-2)}. Let $\gamma \in (0,1]$. If $[u_{0},u_{1}] \in H^{1}({\bf R}^{n})\times (L^{2}({\bf R}^{n})\cap L^{1,\gamma}({\bf R}^{n}))$ further satisfies
\[\int_{{\bf R}^{2}}(1+\vert x\vert)\left(\frac{1}{c(x)^{2}}\vert u_{1}(x)\vert^{2} + \vert\nabla u_{0}(x)\vert^{2}\right)dx < +\infty,\]
and
\[\int_{{\bf R}^{2}}\frac{u_{1}(x)}{c(x)^{2}}dx = 0,\]
then the unique solution $u \in C_{1}^{2}$ to problem {\rm (1.1)}-{\rm (1.2)} satisfies 
\[E_{R}(t) = O(t^{-(1-\eta)})\quad (t \to \infty),\]
for each $R > L$ provided that $\eta := 2L\Vert\frac{1}{c(\cdot)}\Vert_{\infty}\Vert\nabla c\Vert_{\infty} \in [0,1)$.
\end{theo}
\begin{theo}\,Let $n = 1$, and assume {\rm (A-1)} and {\rm (A-2)}. If $[u_{0},u_{1}] \in H^{1}({\bf R}^{n})\times L^{2}({\bf R}^{n})$ further satisfies
\[\int_{{\bf R}}(1+\vert x\vert)\left(\frac{1}{c(x)^{2}}\vert u_{1}(x)\vert^{2} + \vert\nabla u_{0}(x)\vert^{2}\right)dx < +\infty,\]
then the unique solution $u \in C_{1}^{1}$ to problem {\rm (1.1)}-{\rm (1.2)} satisfies 
\[E_{R}(t) = O(t^{-(1-\eta)})\quad (t \to \infty),\]
for each $R > L$ provided that $\eta := 2L\Vert\frac{1}{c(\cdot)}\Vert_{\infty}\Vert\nabla c\Vert_{\infty} \in [0,1)$.
\end{theo}

\begin{rem}{\rm Our gain is that the $n = 1$ dimensional case is included in our results, and we do not assume any compactness of the supports of initial data, and weaker regularity assumptions such as $[u_{0},u_{1}] \in H^{1}({\bf R}^{n})\times L^{2}({\bf R}^{n})$ are imposed. Instead of stronger regularity as in \cite{S} we have to pay a price to assume various weighted conditions on the initial data in some functional spaces, and the parameter $\eta$ must be chosen to satisfy $\eta \in [0,1)$. This condition on $\eta$ is crucial in this paper. In this connection, if $c(x) = 1$ for all $x \in {\bf R}^{n}$, then $\Vert\nabla c\Vert_{\infty} = 0$, so that $\eta = 0$, and in this case the obtained results remind us of those of \cite{IN} studied in an exterior domain with a star-shaped compliment set  (star-shaped obstacle).
}  
\end{rem}

\begin{rem}{\rm For example, if $L > 0$ is small, $\displaystyle{\inf_{x \in {\bf R}^{n}}}c(x)$ is sufficiently far from $0$, and $\Vert\nabla c\Vert_{\infty}$ is small, then we can realize the hypothesis $\eta \in [0,1)$. The smallness of $L$ implies $-c(x)^{2}\Delta = -\Delta$ for $x \in {\bf R}^{n}\setminus B(0,\varepsilon)$ with small $\varepsilon > 0$. Note that $\displaystyle{\inf_{x \in {\bf R}^{n}}}c(x) > 0$ under the assumption (A-1).}
\end{rem}

\par
\vspace{0.2cm}

This paper is organized as follows. In section 2 after preparing several propositions and lemmas we shall prove Theorems 1.2, 1.3 and 1.4 at a stroke. The key tool is already prepared in \cite{IS}.\\

{\bf Notation.} {\small Throughout this paper, $\| \cdot\|_q$ stands for the usual $L^q({\bf R}^{n})$-norm. For simplicity of notation, in particular, we use $\| \cdot\|$ instead of $\| \cdot\|_2$. Furthermore, we denote $\Vert\cdot\Vert_{H^{1}}$ as the usual $H^{1}$-norm. On the other hand, we denote the Fourier transform of $f(x)$ by 
\[{\cal F}(f)(\xi) := \hat{f}(\xi) := (\displaystyle{\frac{1}{2\pi}})^{\frac{n}{2}}\displaystyle{\int_{{\bf R}^{n}}}e^{-ix\cdot\xi}f(x)dx\]
as usual with $i := \sqrt{-1}$. As the $L^{2}$-inner product, one employs the following notation:
\[(f,g) := \int_{{\bf R}^{n}}f(x)g(x)dx, \quad f,g \in L^{2}({\bf R}^{n}).\] 
}

\section{$L^2$-bounds of solutions}

In order to prove the previous theorems one first prepare in this section the so-called Morawetz identity. This is our starting point.
\begin{pro}\,Let $n \geq 1$. Under the assumption {\rm (A-1)}, the (unique) weak solution $u \in C_{1}^{n}$ to problem {\rm (1.1)}-{\rm (1.2)} satisfies
\[tE_{u}(t) = \frac{n-1}{2}(\frac{1}{c(\cdot)^{2}}u_{1}, u_{0}) + (\frac{1}{c(\cdot)^{2}}u_{1}, x\cdot\nabla u_{0})\] 
\[-\frac{n-1}{2}(\frac{1}{c(\cdot)^{2}}u_{t}(t,\cdot), u(t,\cdot)) - (\frac{1}{c(\cdot)^{2}}u_{t}(t,\cdot), x\cdot\nabla u(t,\cdot))\]
\[ + \int_{0}^{t}\int_{{\bf R}^{n}}\frac{1}{c(x)^{3}}(x\cdot\nabla c(x))\vert u_{s}(s,x)\vert^{2}dxds\quad (t \geq 0).\]
\end{pro}
The proof of the Morawetz identity can be derived first to the smooth solution  with $u(t,x)$ for  initial data with compact support, say $[u_{0}, u_{1}] \in C_{0}^{\infty}({\bf R}^{n})\times C_{0}^{\infty}({\bf R}^{n})$, by relying on  the multiplier
$$M(u):=tu_t +  x\cdot \nabla u + \frac{n-1}{2}u,$$
 the finite speed of propagation property,  integration by parts, and then by the density arguments. The final identity can be established to the desired weak solution $u \in C_{1}^{n}$. Note that $c \in C({\bf R}^{n}) \cap L^{\infty}({\bf R}^{n})$ under the assumption (A-1) (cf. \cite[Theorem IX.12]{HB}).

As a second work, we derive several $L^{2}$-bounds of solutions under non-compact support conditions on the initial data. For this purpose, we rely on an improvement version of an original idea established in \cite{IM}  because we can now use the Fourier transform appropriately to obtain them. We have the following significant
propositions.  These results will  be used when one estimates the term $(\frac{1}{c(\cdot)^{2}}u_{t}(t,\cdot), u(t,\cdot))$ in Proposition 2.1.

\begin{pro} Let $n \geq 3$. If $[u_{0},u_{1}] \in H^{1}({\bf R}^{n})\times (L^{2}({\bf R}^{n})\cap L^{1}({\bf R}^{n}))$, then the unique solution $u \in C_{1}^{n}$ to problem {\rm (1.1)}-{\rm (1.2)} satisfies
\[\Vert u(t,\cdot)\Vert \leq C\Vert  c^{-1}\Vert_{\infty}^{2}\big(\Vert u_{1}\Vert + \Vert u_{ 1}\Vert_{1}\big) + C\Vert c^{-1} u_{0}\Vert,\]
with some constant $C > 0$.
\end{pro}

\begin{pro} Let $n = 2$ and $\gamma \in (0,1]$. If $[u_{0},u_{1}] \in H^{1}({\bf R}^{n})\times (L^{2}({\bf R}^{n})\cap L^{1,\gamma}({\bf R}^{n}))$ further satisfies 
\[\int_{{\bf R}^{n}}\frac{u_{1}(x)}{c(x)^{2}}dx = 0,\]
then the unique solution $u \in C_{1}^{2}$ to problem {\rm (1.1)}-{\rm (1.2)} satisfies
\[\Vert u(t,\cdot)\Vert \leq C\Vert c^{-1}\Vert_{\infty}^{2}\big(\Vert u_{1}\Vert + \Vert u_{1}\Vert_{1,\gamma}\big) + C\Vert \bf c^{-1}u_{0}\Vert,\]
with some constant $C > 0$.
\end{pro}

In the course of proofs of Propositions 2.2 and 2.3, the next inequality concerning the Fourier image of the Riesz potential plays an crucial role. This comes from \cite[(ii) of Proposition 2.1]{Ike-0}. 
\begin{pro}\,Let $[n,\gamma,\theta]$ satisfy $n \geq 1$, $\gamma \in [0,1]$ and $\theta \in [0,\gamma + \displaystyle{\frac{n}{2}})$. Then, for all $f \in L^{2}({\bf R}^{n})\cap L^{1,\gamma}({\bf R}^{n})$ satisfying
\[\int_{{\bf R}^{n}}f(x)dx = 0\]
it is true that 
\[\displaystyle{\int_{{\bf R}^{n}}}\displaystyle{\frac{\vert \hat{f}(\xi)\vert^{2}}{\vert\xi\vert^{2\theta}}}d\xi \leq C(\Vert f\Vert_{1,\gamma}^{2} + \Vert f\Vert^{2}) \]
with some constant $C = C_{n,\theta,\gamma} > 0$.
\end{pro}

While, the following result is well-known (cf., B. Muckenhoupt \cite[Theorem 1]{Mu})
\begin{pro}\,Let $[n,\gamma,\theta]$ satisfy $n \geq 1$, $\gamma \in [0,1]$ and $\theta \in [0,\displaystyle{\frac{n}{2}})$. Then, for all $f \in L^{2}({\bf R}^{n})\cap L^{1,\gamma}({\bf R}^{n})$ it is true that 
\[\displaystyle{\int_{{\bf R}^{n}}}\displaystyle{\frac{\vert \hat{f}(\xi)\vert^{2}}{\vert\xi\vert^{2\theta}}}d\xi \leq C\big(\;\Vert f\Vert_{1,\gamma}^{2} + \Vert f\Vert^{2} + \left\vert\int_{{\bf R}^{n}}f(x)dx\right\vert^{2}\big) \]
with some constant $C = C_{n,\theta,\gamma} > 0$.
\end{pro}
{\it \bf Proof of Propositions 2.2 and 2.3.}\, Let us prove Propositions 2.2 and 2.3 at a stroke. We use the idea from \cite{IM} in the Fourier space. Let $v(t,x) := \displaystyle{\int_{0}^{t}}u(s,x)ds$. Then, the function $v$ satisfies
\begin{equation}
\frac{1}{c(x)^{2}}v_{tt}(t,x) -\Delta v(t,x) = \frac{1}{c(x)^{2}}u_{1}(x),
\end{equation}
\begin{equation}
v(0,x)= 0,\ \ v_{t}(0,x)= u_{0}(x),\ \ \ x\in{\bf R}^{n}.
\end{equation}
Multiplying both sides of (2.1) by $v_{t}$, and integrating over $[0,t]$  we  derive that
\begin{equation}\label{Ev}
E_{v}(t) = \frac{1}{2}\Vert  c^{-1} u_{0}\Vert^{2} + \int_{{\bf R}^{n}}w(x)v(t,x)dx.
\end{equation} 
where
\[w(x) := \frac{u_{1}(x)}{c(x)^{2}}.\]
Note that $w \in L^{1}({\bf R}^{n})\cap L^{2}({\bf R}^{n})$ because of the assumption (A-1)  and the condition that $u_1 \in L^1 \cap L^2$. By using the Plancherel theorem, the last term of (2.3) can be estimated as follows. For any $\varepsilon > 0$ there exists a constant $C_{\varepsilon} > 0$ such that for all $t \geq 0$ one has
\[\left\vert\int_{{\bf R}^{n}}w(x)v(t,x)dx\right\vert = \left\vert \int_{{\bf R}_{\xi}^{n}}\hat{w}(\xi)\overline{\hat{v}(t,\xi)}d\xi\right\vert\]
\[\leq \int_{{\bf R}_{\xi}^{n}}\vert\hat{w}(\xi)\vert\vert\hat{v}(t,\xi)\vert d\xi\]
\[\leq C_{\varepsilon}\int_{{\bf R}_{\xi}^{n}}\frac{\vert\hat{w}(\xi)\vert^{2}}{\vert\xi\vert^{2}}d\xi + \varepsilon\Vert\nabla v(t,\cdot)\Vert^{2}.\]

Thus, one has from \eqref{Ev}
\begin{equation}
\frac{1}{2}\Vert v_{t}(t,\cdot)\Vert^{2} + \big(\frac{1}{2}-\varepsilon\big) \Vert\nabla v(t,\cdot)\Vert^{2} \leq \frac{1}{2}\Vert \bf c^{-1}u_{0}\Vert^{2} + C_{\varepsilon}\int_{{\bf R}_{\xi}^{n}}\frac{\vert\hat{w}(\xi)\vert^{2}}{\vert\xi\vert^{2}}d\xi.
\end{equation}

Now, when $n \geq 3$, by using Proposition 2.5 with $\gamma = 0$, one has
\[\int_{{\bf R}_{\xi}^{n}}\frac{\vert\hat{w}(\xi)\vert^{2}}{\vert\xi\vert^{2}}d\xi \leq C\big(\;\Vert w\Vert_{1}^{2} + \Vert w\Vert^{2}\big)\]
\begin{equation}
\leq C\Vert\frac{1}{c} \Vert_{\infty}^{4}\big(\;\Vert u_{1}\Vert_{1}^{2} + \Vert u_{1}\Vert^{2}\big).
\end{equation}
In the case when $n = 2$, by relying on Proposition 2.4 with $\gamma \in (0,1]$ for $n = 2$, one can have 
\begin{equation}
\int_{{\bf R}_{\\xi}^{n}}\frac{\vert\hat{w}(\xi)\vert^{2}}{\vert\xi\vert^{2}}d\xi \leq C\Vert\frac{1}{c} \Vert_{\infty}^{4}(\Vert u_{1}\Vert_{1,\gamma}^{2} + \Vert u_{1}\Vert^{2}),
\end{equation}
where one has just used the assumption $\displaystyle{\int_{{\bf R}^{n}}}w(x)dx = 0$  in Proposition 2.3. Therefore, it follows from (2.4), (2.5) or (2.6) that
\[\frac{1}{2}\Vert v_{t}(t,\cdot)\Vert^{2} + (1-\varepsilon)\Vert\nabla v(t,\cdot)\Vert^{2}\]
\[\leq \frac{1}{2}\Vert \textcolor{blue}{\bf c^{-1}} u_{0}\Vert^{2} + C\Vert c^{-1}\Vert_{\infty}^{4}(\Vert u_{1}\Vert_{1}^{2} + \Vert u_{1}\Vert^{2}), \quad (n \geq 3),\]
and
\[\frac{1}{2}\Vert v_{t}(t,\cdot)\Vert^{2} + (1-\varepsilon)\Vert\nabla v(t,\cdot)\Vert^{2}\]
\[\leq \frac{1}{2}\Vert \textcolor{blue}{\bf c^{-1}}u_{0}\Vert^{2} + C\Vert c^{-1} \Vert_{\infty}^{4}(\Vert u_{1}\Vert_{1,\gamma}^{2} + \Vert u_{1}\Vert^{2}),\quad (n = 2, \gamma \in (0,1]).\]
The desired estimates  of the Propositions 2.2 and 2.3 can be derived because of $v_{t} = u$ by taking $\varepsilon > 0$ small enough.
\hfill
$\Box$
\par 
\vspace{0.2cm}

 The next step  is to handle with the term $(\frac{1}{c(\cdot)^{2}}u_{t}(t,\cdot), x\cdot\nabla u(t,\cdot))$ which appears in Morawetz identity. To do that we prepare an important weighted energy estimate below. For this purpose we define  a weight function $\psi:[0,\infty)\times{\bf R}^{n} \to {\bf R}$  given by
\[\psi(t,x) = \left\{
  \begin{array}{ll}
  \displaystyle{(1+\vert x\vert-t)}&
      \qquad (\vert x\vert \geq t), \\[0.2cm]
  \displaystyle{(1+t-\vert x\vert)^{-1}}& \qquad (\vert x\vert < t).
 \end{array} \right. \]
The following properties of the function $\psi \in C^{1}([0,\infty)\times{\bf R}^{n})$ are a direct calculation
\begin{equation}
\frac{\partial\psi}{\partial t}(t,x) < 0, \quad (t,x) \in [0,\infty)\times{\bf R}^{n}, 
\end{equation} 
\begin{equation}\label{Eiko}
\vert\nabla \psi(t,x)\vert^{2} - (\frac{\partial\psi}{\partial t}(t,x))^{2} = 0, \quad (t,x) \in [0,\infty)\times{\bf R}^{n}, 
\end{equation}
and 
\begin{equation}
\psi(t,x) > 0, \quad (t,x) \in [0,\infty)\times{\bf R}^{n}. 
\end{equation}
The equation \eqref{Eiko} is the so-called Eikonal equation to the free wave equation
\[u_{tt}(t,x) -\Delta u(t,x) = 0.\]
How to choose $\psi(t,x)$ has its origin in \cite{IN}. Based on the (modified) weighted energy estimates originated in \cite{TY}, one can get the following lemma.

\begin{lem}\,Let $n \geq 1$ and assume {\rm (A-2)}. If the initial data $[u_{0},u_{1}] \in H^{1}({\bf R}^{n}) \times L^{2}({\bf R}^{n})$ further satisfies
\begin{equation}
I_{0}^{2} := \int_{{\bf R}^{n}}(1+\vert x\vert)\left(\frac{1}{c(x)^{2}}\vert u_{1}(x)\vert^{2} + \vert\nabla u_{0}(x)\vert^{2}\right)dx < +\infty,
\end{equation}
for each $R > L$, then  the solution $u \in C_{1}^{n}$ to the  problem {\rm (1.1)-(1.2)} satisfies
\[\int_{\vert x\vert \geq R}\psi(t,x)\left(\frac{1}{c(x)^{2}}\vert u_{t}(t,x)\vert^{2} + \vert\nabla u(t,x)\vert^{2}\right)dx \leq CI_{0}^{2}\quad (t \geq 0),\]
where the constant  $C=2+ L > 0$ .
\end{lem}  

Since one can not apply the style itself of the weighted function $\psi(t,x)$ to the equation (1.1) with variable coefficient $c(x)$, one has to introduce a new auxiliary function $\phi(t)$ of $\psi(t,x)$  defined by
\[\phi(t) = \left\{
  \begin{array}{ll}
  \displaystyle{(1+L-t)}&
      \qquad (L > t), \\[0.2cm]
  \displaystyle{(1+t-L)^{-1}}& \qquad (L \leq t).
 \end{array} \right. \]  
Note that $\phi \in C^{1}([0,\infty))$, and
\begin{equation}
\phi_{t}(t) < 0.
\end{equation}
It is easy to check that $\psi(t,x)$ and $\phi(t)$ satisfy the following identities.
\[0 = (\psi u_{t})\left(\frac{1}{c(x)^{2}}u_{tt}-\Delta u\right) = \frac{d}{dt}(\psi(t,x)E(t,x)) - \nabla\cdot(\psi u_{t}\nabla u)\]
\begin{equation}
-\frac{1}{2\psi_{t}}\vert\psi_{t}\nabla u-u_{t}\nabla\psi\vert^{2} + \frac{u_{t}^{2}}{2c(x)^{2}\psi_{t}}\left(c(x)^{2}\vert\nabla\psi\vert^{2}-\psi_{t}^{2}\right),
\end{equation}
and
\[0 = (\phi u_{t})\left(\frac{1}{c(x)^{2}}u_{tt}-\Delta u\right) = \frac{d}{dt}(\phi(t)E(t,x)) - \phi(t)\nabla\cdot(u_{t}\nabla u)\]
\begin{equation}
- \frac{\phi_{t}(t)}{2c(x)^{2}}\left(c(x)^{2}\vert\nabla u\vert^{2} +u_{t}^{2}\right),
\end{equation}
where 
\begin{equation}
E(t,x) := \frac{1}{2}\left(\frac{1}{c(x)^{2}}\vert u_{t}(t,x)\vert^{2} + \vert\nabla u(t,x)\vert^{2} \right).
\end{equation}

{\it \bf   Proof of Lemma 2.1.} The proof  can be proceeded similarly to \cite[Lemma 3.3]{IS} with a constant $r_{0}$ replaced by $L$. However, for making this paper self-contained we shall draw its full proof.

\,It follows from (2.7) and (2.12) that
\[0 \geq \frac{d}{dt}(\psi(t,x)E(t,x)) - \nabla\cdot(\psi u_{t}\nabla u)\]
\begin{equation}
+ \frac{u_{t}^{2}}{2c(x)^{2}\psi_{t}}\left(c(x)^{2}\vert\nabla\psi\vert^{2}-\psi_{t}^{2}\right).
\end{equation}
Integrating (2.15) over $[0,t]\times ({\bf R}^{n}\setminus B(0,L))$ one can get
\[\int_{\vert x\vert \geq L}\psi(0,x)E(0,x)dx + \int_{0}^{t}\int_{\vert x\vert \geq L}\nabla\cdot(\psi(s,x)u_{s}(s,x)\nabla u(s,x))dxds\]
\[\geq \int_{\vert x\vert \geq L}\psi(t,x)E(t,x)dx +  \int_{0}^{t}\int_{\vert x\vert \geq L}\frac{u_{s}(s,x)^{2}}{2\psi_{s}(s,x)}(\vert\nabla\psi(s,x)\vert^{2}-\psi_{s}(s,x)^{2})dxds,\]
where one has just used the assumption (A-2). By applying (2.8) one obtains
\[\int_{\vert x\vert \geq L}\psi(0,x)E(0,x)dx + \int_{0}^{t}\int_{\vert x\vert \geq L}\nabla\cdot(\psi(s,x)u_{s}(s,x)\nabla u(s,x))dxds\]
\begin{equation}
\geq \int_{\vert x\vert \geq L}\psi(t,x)E(t,x)dx.
\end{equation}
On the other hand, by integrating (2.13) over $[0,t]\times B(0,L)$, because of (2.11) one can get
\[\int_{\vert x\vert \leq L}\phi(0)E(0,x)dx + \int_{0}^{t}\int_{\vert x\vert \leq L}\phi(s)\nabla\cdot(u_{s}(s,x)\nabla u(s,x))dxds\]
\begin{equation}
\geq \int_{\vert x\vert \leq L}\phi(t)E(t,x)dx.
\end{equation}
Now, since $\phi(t) = \psi(t,x)$ on the sphere $\vert x\vert = L$, it follows from the divergence formula that
\[\int_{0}^{t}\int_{\vert x\vert \geq L}\nabla\cdot(\psi(s,x)u_{s}(s,x)\nabla u(s,x))dxds + \int_{0}^{t}\int_{\vert x\vert \leq L}\phi(s)\nabla\cdot(u_{s}(s,x)\nabla u(s,x))dxds = 0.\]
Thus, by summing up (2.16) and (2.17) one can arrived at the desired estimate
\[\int_{\vert x\vert \geq R}\psi(t,x)E(t,x)dx \leq \int_{\vert x\vert \geq L}(1+\vert x\vert)E(0,x)dx + (1+L)\int_{\vert x\vert \leq L}E(0,x)dx,\]
where one has used the fact that $R > L$ and $\phi(t) > 0$.
\hfill
$\Box$
\par 
\vspace{0.2cm}

Based on Lemma 2.1, one can also obtain the following lemma. This is re-stated version of \cite[Lemma 3.4]{IS} with $r_{0}$ replaced by $L$. Note that from the assumption on $c(x)$, one has
\[c_{m} := \inf_{x \in {\bf R}^{n}}c(x) > 0.\]

\begin{lem} Let $R > L$,  $t > R$, and $c(x)$ satisfies the assumptions {\rm (A-1)} and {\rm (A-2)}. Then  it is true that
\[\left\vert\left(\frac{1}{c(\cdot)^{2}}u_{t}(t,\cdot), x\cdot\nabla u(t,\cdot)\right)\right\vert
\leq \frac{R}{c_{m}}E_{R}(t) + C\frac{I_{0}^{2}}{2} + t\int_{\vert x\vert \geq R}E(t,x)dx,\]
where $C > 0$ is a constant independent from initial data.
\end{lem} 
\begin{rem}{\rm When one checks the proof of Lemma 2.2 above, one has to use the assumption (A-2) such that $c(x) = 1$ for $x \in {\bf R}^{n}$ satisfying $\vert x\vert > L$, and in this case one notices that (see (2.14))
\[E(t,x) = \frac{1}{2}\left(\vert u_{t}(t,x)\vert^{2} + \vert\nabla u(t,x)\vert^{2}\right),\quad \vert x\vert > L, \quad t \geq 0.\]
}
\end{rem}


\section{Proof of Theorems 1.2, 1.3 and 1.4}
Now, in this section,  we are in a position to prove Theorems 1.2, 1.3 and 1.4 based on Lemmas 2.1 and 2.2.

Note that it suffices to check only the case for $n \geq 3$, and the case for $n = 2$ is similar by using Proposition 2.3 in place of Proposition 2.2. Note also that for the case of $n = 1$, one does not need to use Proposition 2.3 because of the existence of the coefficient $\displaystyle{\frac{n-1}{2}}$ in Proposition 2.1.

Let $R > L$. We first start with Proposition 2.1 under the assumption (A-2):
\[tE_{u}(t) = J_{0}^{2} -\frac{n-1}{2}(\frac{1}{c(\cdot)^{2}}u_{t}(t,\cdot), u(t,\cdot)) - (\frac{1}{c(\cdot)^{2}}u_{t}(t,\cdot), x\cdot\nabla u(t,\cdot))\]
\begin{equation}
+ \int_{0}^{t}\int_{\vert x\vert \leq L}\frac{1}{c(x)^{3}}(x\cdot\nabla c(x))\vert u_{s}(s,x)\vert^{2}dxds\quad (t \geq 0),
\end{equation}
where
\[J_{0}^{2} := \frac{n-1}{2}(\frac{1}{c(\cdot)^{2}}u_{1}, u_{0}) + (\frac{1}{c(\cdot)^{2}}u_{1}, x\cdot\nabla u_{0}).\]
We observe that  $J_0^2$ is not necessarily  positive.

Note that under the assumption on the regularity imposed on the initial data, one can check that (see (2.10))
\[\vert(\frac{1}{c(\cdot)^{2}}u_{1}, x\cdot\nabla u_{0})\vert \leq \Vert\frac{1}{c(\cdot)}\Vert_{\infty}\int_{{\bf R}^{n}}(\frac{1}{c(x)}\sqrt{\vert x\vert}\vert u_{1}(x)\vert)(\sqrt{\vert x\vert}\vert\nabla u_{0}(x)\vert)\]
\[\leq \frac{1}{2}\Vert\frac{1}{c(\cdot)}\Vert_{\infty}\int_{{\bf R}^{n}}\left(\frac{1}{c(x)^{2}}\vert x\vert\vert u_{1}(x)\vert^{2} + \vert x\vert\vert\nabla u_{0}(x)\vert^{2}\right)dx < +\infty.\]
Then, it follows from the Schwarz inequality, (A-1) and (A-2) that
\[\int_{0}^{t}\int_{\vert x\vert \leq L}\frac{1}{c(x)^{3}}(x\cdot\nabla c(x))\vert u_{s}(s,x)\vert^{2}dxds \leq 2L\Vert\frac{1}{c(\cdot)}\Vert_{\infty}\Vert\nabla c\Vert_{\infty}\int_{0}^{t}\int_{\vert x\vert \leq R}\frac{1}{2}\frac{1}{c(x)^{2}}\vert u_{s}(s,x)\vert^{2}dxds\]
\begin{equation}
\leq  2L\Vert\frac{1}{c(\cdot)}\Vert_{\infty}\Vert\nabla c\Vert_{\infty}\int_{0}^{t}E_{R}(s)ds.
\end{equation}
Since
\[tE_{u}(t) = tE_{R}(t) + t\int_{\vert x\vert \geq R}E(t,x)dx,\]
it follows from (3.1) and (3.2) that
\[tE_{R}(t) + t\int_{\vert x\vert \geq R}E(t,x)dx \leq J_{0}^{2} + \frac{n-1}{2}\left\vert (\frac{1}{c(\cdot)^{2}}u_{t}(t,\cdot), u(t,\cdot))\right\vert\]
\begin{equation}
+ \left\vert(\frac{1}{c(\cdot)^{2}}u_{t}(t,\cdot), x\cdot\nabla u(t,\cdot))\right\vert + \eta\int_{0}^{t}E_{R}(s)ds,
\end{equation}
where
\[\eta := 2L\Vert\frac{1}{c(\cdot)}\Vert_{\infty}\Vert\nabla c\Vert_{\infty}.\]
This parameter $\eta$ is quite important.
 
While, by Lemma 2.2 and (3.3) one can obtain
\[tE_{R}(t) + t\int_{\vert x\vert \geq R}E(t,x)dx \leq J_{0}^{2} + \frac{n-1}{2}\left\vert (\frac{1}{c(\cdot)^{2}}u_{t}(t,\cdot), u(t,\cdot))\right\vert\]
\begin{equation}
+ \frac{R}{c_{m}}E_{R}(t) + C\frac{I_{0}^{2}}{2} + t\int_{\vert x\vert \geq R}E(t,x)dx + \eta\int_{0}^{t}E_{R}(s)ds,
\end{equation}
which implies
\begin{equation}
(t-\frac{R}{c_{m}})E_{R}(t) \leq J_{0}^{2} + C\frac{I_{0}^{2}}{2}+ \frac{n-1}{2}\left\vert (\frac{1}{c(\cdot)^{2}}u_{t}(t,\cdot), u(t,\cdot))\right\vert + \eta\int_{0}^{t}E_{R}(s)ds.
\end{equation}
Here, let us estimate the term 
\[\left\vert (\frac{1}{c(\cdot)^{2}}u_{t}(t,\cdot), u(t,\cdot))\right\vert,\]
by relying on Propositions 2.2 and 2.3.

We estimate only the case for $n \geq 3$. Based on Proposition 2.2, (A-1) and the Schwarz inequality  we can derive 
\[\left\vert (\frac{1}{c(\cdot)^{2}}u_{t}(t,\cdot), u(t,\cdot))\right\vert \leq c_{m}^{-1}\Vert\frac{1}{c(\cdot)} u_{t}(t,\cdot)\Vert\Vert u(t,\cdot)\Vert \]
\[\leq \frac{1}{2c_{m}^{2}}\Vert\frac{1}{c(\cdot)} u_{t}(t,\cdot)\Vert^{2} + \frac{1}{2}\Vert u(t,\cdot)\Vert^{2}\]
\[\leq \frac{1}{2c_{m}^{2}}\Vert\frac{1}{c(\cdot)} u_{t}(t,\cdot)\Vert^{2} + C\Vert\frac{1}{c}\Vert_{\infty}^{4}(\Vert u_{1}\Vert + \Vert u_{1}\Vert_{1})^{2} + C\Vert u_{0}\Vert^{2}.\]
Since $2^{-1}\Vert\frac{1}{c(\cdot)}u_{t}(t,\cdot)\Vert^{2} \leq E_{u}(t) = E_{u}(0)$ (see (1.3)), one has
\begin{equation}
\left\vert (\frac{1}{c(\cdot)^{2}}u_{t}(t,\cdot), u(t,\cdot))\right\vert \leq \frac{1}{c_{m}^{2}}E_{u}(0) + C\Vert\frac{1}{c}\Vert_{\infty}^{4}(\Vert u_{1}\Vert + \Vert u_{1}\Vert_{1})^{2} + C\Vert u_{0}\Vert^{2}.
\end{equation}

Because of (3.5) and (3.6) one can arrive at the significant inequality of the Gronwall type:
\[(t-\frac{R}{c_{m}})E_{R}(t) \leq J_{0}^{2} + C\frac{I_{0}^{2}}{2}+ \frac{n-1}{2}\left\vert (\frac{1}{c(\cdot)^{2}}u_{t}(t,\cdot), u(t,\cdot))\right\vert\]
\begin{equation}
+ \eta\int_{0}^{t}E_{R}(s)ds \leq K_{0}^{2} + \eta\int_{0}^{t}E_{R}(s)ds,
\end{equation}
where
\[K_{0}^{2} := J_{0}^{2} + C\frac{I_{0}^{2}}{2} + \frac{n-1}{2}\frac{1}{c_{m}^{2}}E_{u}(0) + C\frac{n-1}{2}\Vert\frac{1}{c}\Vert_{\infty}^{4}(\Vert u_{1}\Vert + \Vert u_{1}\Vert_{1})^{2} + C\frac{n-1}{2}\Vert u_{0}\Vert^{2}.\]

Now, let us solve the integral inequality (3.7) under the assumption $\eta \in [0,1)$. This is rather standard. For completeness we write its full proof.

To do this we consider the function
\[\xi(t) := (t-\frac{R}{c_{m}})^{-\eta}\int_{0}^{t}E_{R}(s)ds\quad (t > R/c_{m}).\]

Then, it follows from (3.7) that
\[\xi'(t) = (t-\frac{R}{c_{m}})^{-1-\eta}\left\{\left(t-\frac{R}{c_{m}}\right)E_{R}(t) -\eta\int_{0}^{t}E_{R}(s)ds\right\} 
\leq K_{0}^{2}(t-\frac{R}{c_{m}})^{-1-\eta}, \quad (t > R/c_{m}).\]

Integrating over $[t_{0},t]$ with large $t_{0} \gg 1$, one can get
\begin{equation}\xi(t) \leq \xi(t_{0}) + K_{0}^{2}\int_{t_{0}}^{t}(s-\frac{R}{c_{m}})^{-1-\eta}ds
\leq \;  \xi(t_{0}) + K_{0}^{2}\eta^{-1}\big(\;t_0-\frac{R}{c_{m} }\big)^{-\eta} =: M_{0}.
\end{equation}
By (3.7) and (3.8) one can obtain the desired estimate
\[(t-\frac{R}{c_{m}})E_{R}(t) \leq K_{0}^{2} + \eta\int_{0}^{t}E_{R}(s)ds \leq K_{0}^{2} + \eta M_{0}(t-\frac{R}{c_{m}})^{\eta}\quad (t > t_{0} \gg 1).\]
This completes the proof of Theorems 1.2, 1.3 and 1.4.
\hfill
$\Box$
\par 
\vspace{0.2cm}

\par
\par
\vspace{0.5cm}
\noindent{\em Acknowledgement.}
\smallskip
The work of the second author (R. IKEHATA) was supported in part by Grant-in-Aid for Scientific Research (C)15K04958  of JSPS. 


\end{document}